\documentclass[11pt,reqno,oneside]{amsart}
\usepackage{geometry}                
\usepackage[parfill]{parskip}    
\usepackage{graphicx}
\usepackage{amssymb}
\usepackage{epstopdf}
\usepackage{tikz}

\usepackage[frozencache]{minted}

\newtheorem{prop}{Proposition}[section]

\title{A Census of Small Transitive Groups and Vertex-Transitive Graphs}
\author{Derek Holt}
\address{Mathematics Institute, University of Warwick}
\email{D.F.Holt@warwick.ac.uk}

\author{Gordon Royle}
\address{Mathematics and Statistics, University of Western Australia}
\email{gordon.royle@uwa.edu.au}


\begin{document}

\begin{abstract}
We describe two similar but independently-coded computations used to construct a complete catalogue of the transitive groups of degree less than $48$, thereby verifying, unifying and extending the catalogues previously available. From this list, we construct all the vertex-transitive graphs of order less than $48$.  We then present a variety of summary data regarding the transitive groups and vertex-transitive graphs, focussing on properties that seem to occur most frequently in the study of groups acting on graphs.  We illustrate how such catalogues can be used, first by finding a complete list of the elusive groups of order at most $47$ and then by completely determining which groups of order at most $47$ are CI groups. 
\end{abstract}

\maketitle

\section{Introduction}

In the study of finite permutation groups and the study of groups acting on graphs, transitive groups play a fundamental role.  For well over 100 years, researchers have used catalogues of transitive groups of (necessarily) small degree to 
provide examples and counterexamples, or simply to identify promising lines of enquiry.  
G. A. Miller  was one
of the most prolific of these early cataloguers with many of the 400+ papers in his \emph{Collected Works} \cite{GAMillerCollected} listing small groups with particular properties (transitive, primitive, etc).  The difficulty of the task is determined by the prime 
factorisation of the degree, with the most difficult degrees being those with several small prime factors. The first ``difficult degree'' is $12$ and although Miller \cite{Miller96} published a supposedly complete list of the transitive permutation 
groups of degree $12$ as early as $1896$, it was later discovered to contain a handful of mistakes (Royle \cite{MR922391}).  Miller notwithstanding, constructing catalogues of combinatorial objects is a 
difficult and error-prone task by hand, yet one particularly amenable to computation.  As computers became more
widespread and more powerful, the catalogues of transitive groups were extended to larger and larger degrees 
until by 1996 the lists were complete up to degree $31$ (see Hulpke \cite{MR2168238} for an overview including references to the various authors who contributed to these extensions).

When the degree is a power of a small prime, then we expect the computation to be the most difficult and to yield the largest collections of transitive groups, and so $32$ is the next difficult case. This case was resolved by Cannon \& Holt \cite{MR2455702} who determined that there are $2801324$ transitive groups of degree $32$, which is two 
orders of magnitude more than the total number of all transitive groups of degree less than $32$. 
(Throughout this paper, when we talk about classifying transitive permutation
groups of a given degree $n$, then we always mean classifying them up to
conjugacy in $\text{Sym}(n)$, and when we say, for example, that there are 
$2801324$ transitive groups of degree $32$, we really mean that there are
$2801324$ conjugacy classes of such subgroups in  $\text{Sym}(n)$. It would
not be difficult to count the sizes of each of these classes and thereby count
the total number of groups.)

In this paper we extend the catalogues of transitive groups up to degree $47$, i.e. just short of the next really difficult degree, namely $48$. For ``end-users'' of a catalogue of combinatorial objects, the primary concern is the extent to which the user has confidence that the catalogue is correct and complete. Increasingly, computational papers in combinatorics are addressing these concerns through redundancy --- two or more researchers writing programs that are independent (as far as possible) and then 
cross-checking the results. Along these lines, we present two computations that were conducted totally independently, and yielded the same output when first compared.
(More precisely, we checked that the lists of groups output by the two
computations had the same lengths, and that each group on one of the lists
was conjugate in $\text{Sym}(n)$ to exactly one group on the other list.)
Both computations were done in {\sc Magma \cite{MR1484478}} and so both rely on basic functions such as calculating the maximal subgroups of a finite group and checking whether two groups are conjugate in the symmetric group. As these functions have by now received decades of user-testing, we feel confident that the chance of an error from this source is negligible. 


The first author is in the process of computing the transitive groups of
degree $48$, but with no great urgency, so this computation is not
expected to complete in the foreseeable future.  Based on the numbers of groups
found so far, it is looking as though there could be about $200$ million such
groups in total, of which all except about $3.4$ million will have minimal
block size $2$. 

\section{First Computation}

The first computation was based around recursive use of the  {\tt MaximalSubgroups} command of {\sc Magma} \cite{MR1484478}, together with partitioning the computation into a large number of disjoint parts that can be run independently (for example, on multiple computers or just multiple cores on a single computer.) It uses only elementary permutation group theory, instead seeking to be as simple and robust as possible. The majority of the effort in 
this approach involves implementing computational techniques to partition the search, to eliminate as much
duplication as possible, and to manage and collate the resulting large data sets; it involves more ``data bookkeeping'' than mathematics.

The \emph{primitive} permutation groups are known up to degree $4095$ (See \cite{MR2845584}) and are incorporated
into the databases of both {\sc Magma} and GAP, and so we need only consider the imprimitive groups. 
By definition if a group $G$ acting transitively on the set $\Omega$ of size $n$ is not primitive, 
then there is at least one partition of $\Omega$ into a \emph{block system} $\mathcal{B}$ such that 
$G$ permutes the blocks of $\mathcal{B}$. If we let $G^\mathcal{B}$ denote the
action of $G$ on the blocks (``the top group'') and if $\mathcal{B}$ has $n/k$ blocks of size $k$, 
then $G^\mathcal{B}$ is a transitive permutation group of degree $n/k$. 
Thus we can associate to each group $G$ a set of pairs of the form
\[
\{(k, G^\mathcal{B}) : \mathcal{B}\text{ is a minimal block system for } G \text{ with blocks of size } k\}.
\]

If this set contains more than one pair (imprimitive groups may of course have more than one minimal block system),
then we wish to distinguish just one of them. Thus we define the \emph{signature} of an imprimitive 
permutation group to be the \emph{lexicographically least}  pair $(k,G^\mathcal{B})$ associated with $G$, 
where the second component is indexed according to its order in the list of transitive groups of 
degree $n/k$ already in {\sc Magma}. 
But note that it can happen that two different minimal block systems of $G$
define the same signature.

Our aim will be to separate the computation into parts, with each part constructing only the groups with a particular signature. Given an integer $k$ such that $1 < k < n$, and a transitive group $H$ of degree $n/k$, the wreath product $\text{Sym}(k) \wr H$ contains (a conjugate of) every transitive group of degree $n$ with signature $(k,H)$. Therefore, at least in principle, these groups can all be found by exploring the subgroup lattice of $\text{Sym}(k) \wr H$. To do this, we repeatedly use the {\tt MaximalSubgroups} command of {\sc Magma} thereby traversing the subgroup lattice downwards and in a breadth-first fashion, pruning each branch of the search as soon as it produces groups with signature differing from $H$, while using conjugacy tests to avoid duplication as far as possible.


More precisely, the program maintains a list of \emph{active} transitive subgroups of $\text{Sym}(k) \wr H$ in descending group order. This list is initialised with the single group $\text{Sym}(k) \wr H$ which by default is 
constructed with the ``canonical block system'' with blocks $\{1,2,\ldots,k\}$, $\{k+1, k+2,\ldots, 2k\}$, and so on. At each step, the program removes an active subgroup of largest order, say $G$, from the list for processing. The processing stage computes the transitive maximal subgroups of $G$, and tests each of these maximal subgroups $M$ to determine if it should be retained or not, using the following rules:
\begin{enumerate}
\item If $M$ has a minimal block system with blocks of size smaller than $k$, then reject,
\item If the action of $M$ on the canonical block system is not equal to $H$, then reject,
\item If $M$ is conjugate (in $\text{Sym}(k) \wr \text{Sym}(n/k)$) to a group on the list of active subgroups,
then reject.
\end{enumerate}
If $M$ passes all these tests, it is then added to the list of active subgroups in the appropriate position (depending on its order).

The rationale behind the first two tests is that if $M$ fails either of them, then a conjugate of $M$ will be found during another part of the search, either the search for all transitive groups preserving a block system with strictly smaller blocks, or during the search for transitive groups with a different ``top group''. The third test ensures that the groups in each 
part are pairwise non-conjugate.

Every group that is processed is a transitive group with a minimal block system with blocks of size $k$, no block systems with blocks with fewer than $k$ elements, and with top-group isomorphic to $H$.  Conversely, every transitive group with these three properties will at some stage enter the list of active subgroups and then be processed. Therefore by adding an output step as the point that a group enters (or leaves) the list of active subgroups, we can compute all transitive groups associated with a particular $(k,H)$ pair. By running over all possibilities for $k$ and letting $H$ range over all the transitive groups of degree $n/k$, every transitive group of degree $n$ is constructed.
As indicated previously, we need to eliminate conjugate groups that occur in
more than one part of the computation, which results from groups having
multiple minimal block systems. If a group has two minimal block
systems that define different signatures, then we only keep it if its
canonical block system has minimal signature. But if there is more than one
minimal block system with the same minimal signature as the canonical block
system, then a final conjugacy check is required to eliminate possible
duplicates.

\section{Second Computation}
For many of the calculations, we used essentially the same methods as in the
first computation, although these calculations were carried out completely
independently. In fact, for all degrees except for $36$ and $40$,
we were able to complete the calculation in a single run, without any
need for filtering by size of blocks, by repeated application
of the {\sc Magma} commands {\tt MaximalSubgroups} and  {\tt IsConjugate},
starting with $\text{Sym}(n)$. But for some degrees, such as $42$ and $44$,
it was quicker to use  the alternative techniques that we shall now describe
for the smaller block sizes.

In the following descriptions,  we use ATLAS notation for group structures.
In particular, an integer  $k$ in a structure description denotes a
cyclic group of order $k$.  For calculations with $n$ even and minimal
block size $2$, we used a method similar to the one described in detail in
\cite[Section 2.2]{MR2455702}, which was applied to the
enumeration of the transitive groups of degree $32$ with minimal block size $2$.
Suppose that the transitive group $G$ preserves the block system
$\mathcal{B}$ with blocks of size $2$. So $G \le W := 2 \wr \text{Sym}(n/2)$.
Then $\bar{H} := G^\mathcal{B}$ is one of the groups in the known list of
transitive groups of degree $n/2$, and $G \le H := 2 \wr \bar{H}$.
As in the first computation, we calculate those groups with signature
$(2,\bar{H})$ for each individual group $\bar{H}$.

Let $K \cong  2^{n/2}$ be the kernel of the action of $W$ on $\mathcal{B}$.
Then we can regard $K$ as a module for $\bar{H}$ over the field
${\mathbb F}_2$ of order $2$, and $M := G \cap K$ is an
${\mathbb F}_2\bar{H}$-submodule. We can use the
{\sc Magma} commands {\tt GModule} and {\tt Submodules} to find all
such submodules.
In fact, since we are looking for representatives of the conjugacy
classes of transitive subgroups of $W$, we only want one representative 
of the conjugation action of $N := 2 \wr N_{\text{Sym}(n/2)}(\bar{H})$
on the set all ${\mathbb F}_2H$-submodules $M$ of $K$, and we
use the {\sc Magma} command {\tt IsConjugate} to find such representatives.

Now, for each such pair $(\bar{H},M)$, the transitive groups $G$ with
$\bar{H} = G^\mathcal{B}$ and $M = G \cap K$ correspond to complements
of $K/M$ in $H/M$, and the $H$-conjugacy classes of such complements
correspond to elements of the cohomology group $H^1(\bar{H},K/M)$, which
can be computed in {\sc Magma}.

We also need to test these groups $G$ for conjugacy under the action of
$N_N(M)$, which can either be done in straightforward fashion using
{\sc Magma}'s {\tt IsConjugate} function or (with a little more programming
but usually faster in terms of computation) using an induced action on the
cohomology group. Finally, for each $G$ that we find,
we need to find all blocks systems with block size $2$ preserved by $G$,
so that we can eliminate occurrences of groups that are conjugate in
$\text{Sym}(n)$ but arise either for distinct pairs $(\bar{H},M)$ or
more than once for the same pair.
We refer the reader to  \cite[Section 2.2]{MR2455702} for further details.

In the case $n=36$, we used a similar method for blocks systems with
blocks of size $3$. In this case, we have $G \le
\text{Sym}(3) \wr \text{Sym}(12) \cong 3^{12}:(2 \wr \text{Sym}(12)).$
Let $H$ be the projection of $G$ onto the quotient group
$2 \wr \text{Sym}(12)$. Then $H \le \text{Sym}(24)$ and $H$
preserves a block system with blocks of size $2$ and projects onto a transitive
subgroup of degree $12$. We can find the possible groups $H$ using the
method just described for blocks of size $2$ (although there is a minor
complication, because $H$ is not necessarily transitive), and then
find the possible groups $G$ using the same method, but working
with modules over ${\mathbb F}_3$ rather than ${\mathbb F}_2$.

In the cases $n=36$ and $n=40$ with minimal block size $4$, we did a
corresponding $3$-step calculation using
$G \le \text{Sym}(4) \wr \text{Sym}(n/4) \cong
2^{n/2}:3^{n/4}:2^{n/4}:\text{Sym}(n/4)$.

Finally, for $n=40$ with minimal block size $5$, the induced action of the
stabilizer of a block on that block is transitive of degree  $5$, and its
structure is one of $5$, $\hbox{5:2}$, $\hbox{5:4}$,
$\text{Alt}(5)$, or
$\text{Sym}(5)$.
We enumerated those groups in the final two of these cases separately using
the methods of the first computation. For the first three cases we used a
$3$-step calculation using $G \le 5^8:2^8\cdot 2^8:\text{Sym}(n/5)$.

These computations were originally carried out in $2014$ on a number of
different computers
with different specifications, so it is difficult to provide a meaningful 
estimate of the total cpu-time involved, but this was of the order of
150 hours in total. Perhaps suprisingly, the most time consuming
individual computation was of the case $n=36$ with block size $9$, which
took about 47 hours. The calculation requiring the most memory was the case
$n=40$ with block size $2$, which used about $27$GB.

\section{Transitive Groups}

The numbers of transitive groups of each degree (up to conjugacy in
$\text{Sym}(n)$)  are shown in Table~\ref{tab:alltrans}. As expected, the number of transitive groups is primarily dependent on the number of repeated prime factors in the degree. The single degree $n=32$ provides well over $90\%$ of the transitive groups in the entire catalogue.

\def\arraystretch{1.2}
\setlength\belowcaptionskip{25pt}

\begin{table}

\begin{tabular}{|cr|cr|cr|cr|cr|}
\hline
$n$&\multicolumn{1}{c|}{$g(n)$} &$n$&\multicolumn{1}{c|}{$g(n)$} &$n$&\multicolumn{1}{c|}{$g(n)$} &$n$&\multicolumn{1}{c|}{$g(n)$} &$n$&\multicolumn{1}{c|}{$g(n)$}\\
\hline
&&$11$&$8$&$21$&$164$&$31$&$12$&$41$&$10$\\
$2$&$1$&$12$&$301$&$22$&$59$&$32$&$2801324$&$42$&$9491$\\
$3$&$2$&$13$&$9$&$23$&$7$&$33$&$162$&$43$&$10$\\
$4$&$5$&$14$&$63$&$24$&$25000$&$34$&$115$&$44$&$2113$\\
$5$&$5$&$15$&$104$&$25$&$211$&$35$&$407$&$45$&$10923$\\
\hline
$6$&$16$&$16$&$1954$&$26$&$96$&$36$&$121279$&$46$&$56$\\
$7$&$7$&$17$&$10$&$27$&$2392$&$37$&$11$&$47$&$6$\\
$8$&$50$&$18$&$983$&$28$&$1854$&$38$&$76$&&\\
$9$&$34$&$19$&$8$&$29$&$8$&$39$&$306$&&\\
$10$&$45$&$20$&$1117$&$30$&$5712$&$40$&$315842$&&\\
\hline
\end{tabular}
\vspace{2mm}
\caption{Numbers $g(n)$ of transitive groups of degree $n$}
\label{tab:alltrans}
\end{table}

\def\arraystretch{1.2}
\begin{table}
\begin{tabular}{|lr|lr|lr|lr|lr|}
\hline
$n$&\multicolumn{1}{c|}{$m(n)$} &$n$&\multicolumn{1}{c|}{$m(n)$} &$n$&\multicolumn{1}{c|}{$m(n)$}&$n$&\multicolumn{1}{c|}{$m(n)$}&$n$&\multicolumn{1}{c|}{$m(n)$}  \\
\hline
&&$11$&$1$&$21$&$5$&$31$&$1$&$41$&$1$\\
$2$&$1$&$12$&$17$&$22$&$6$&$32$&$12033$&$42$&$84$\\
$3$&$1$&$13$&$1$&$23$&$1$&$33$&$3$&$43$&$1$\\
$4$&$2$&$14$&$6$&$24$&$213$&$34$&$7$&$44$&$148$\\
$5$&$1$&$15$&$4$&$25$&$2$&$35$&$4$&$45$&$41$\\
$6$&$4$&$16$&$75$&$26$&$7$&$36$&$436$&$46$&$4$\\
$7$&$1$&$17$&$1$&$27$&$20$&$37$&$1$&$47$&$1$\\
$8$&$5$&$18$&$23$&$28$&$30$&$38$&$5$&&\\
$9$&$2$&$19$&$1$&$29$&$1$&$39$&$4$&&\\
$10$&$6$&$20$&$47$&$30$&$79$&$40$&$1963$&&\\
\hline
\end{tabular}
\vspace{2mm}
\caption{Numbers $m(n)$ of minimal transitive groups of degree $n$}
\label{tab:mintrans}
\end{table}

For many applications that involve considering all possible transitive actions of a certain degree, it is sufficient to consider only the \emph{minimal transitive groups} i.e., transitive groups with no proper transitive subgroups. (One example of this can be found in the next section, where all vertex-transitive graphs are constructed.) Testing if a transitive graph is minimal can be done by finding all its maximal subgroups and verifying that none are transitive. As most of the groups are not minimal transitive, it proves useful in practice to first construct some random subgroups in an attempt to find a transitive proper subgroup, only undertaking the more expensive step of finding all maximal subgroups if this fails.

The numbers of minimal transitive groups of each order are given in Table~\ref{tab:mintrans}, which shows that the numbers of minimal transitive groups vary in much the same way as the numbers of all transitive groups.


\section{Vertex-transitive Graphs}

If $G$ is a group acting transitively on a set $\Omega$, then it is straightforward to construct all the $G$-invariant graphs or digraphs with vertex set $\Omega$. The orbits of $G$ on $\Omega \times \Omega$ are called the \emph{orbitals} of $G$ and it is immediate that the arc set of a $G$-invariant digraph is a union of these orbitals. If $\mathcal{O}$ is an orbital containing a pair $(a,b)$ with $a \ne b$, then we denote by $\mathcal{O}'$ the orbital containing $(b,a)$. It is possible that $\mathcal{O'}$ is equal to $\mathcal{O}$, in which case $\mathcal{O}$ is called \emph{self-paired}, and otherwise $\mathcal{O}$ and $\mathcal{O'}$ are \emph{paired}. Any subset of the orbitals determines a $G$-invariant digraph, and any subset of the orbitals closed under taking pairs determines a $G$-invariant \emph{graph}. (We view a graph simply as a digraph that happens to have the property that if $(a,b)$ is an arc, then so is $(b,a)$, and thus ``digraphs" includes ``graphs".) The set of $G$-invariant digraphs arising in this way will usually contain isomorphic digraphs, but for the small sizes that we are considering, it is easy to filter these out, thus yielding a complete list of the pairwise non-isomorphic $G$-invariant digraphs. 

The complete list of vertex-transitive digraphs can be computed by repeating this computation for each of the transitive groups of degree $n$, combining the resulting lists, and then performing one final filtering process to remove all but one copy of each digraph. To construct only \emph{graphs}, the process is modified slightly to ensure that the orbitals are included/excluded in pairs in the arc-set of the digraph under construction.

If $H$, $G$ are transitive groups such that $H \leq G$, then the orbitals of $G$ are unions of the orbitals of $H$ and so the set of $H$-invariant digraphs contains the set of $G$-invariant digraphs. Thus it is sufficient to perform the construction only for the minimal transitive groups. The \emph{regular} groups of degree $n$ (i.e., those with degree equal to order) are necessarily minimal transitive, and the digraphs arising from these groups are exactly the \emph{Cayley digraphs} of that order. Any non-Cayley digraphs can \emph{only} be created when the larger groups are processed, though these groups will usually produce Cayley digraphs as well.

The numbers of transitive (resp. Cayley) graphs peak at degree $44$ even though there are far fewer groups of degree $44$ than, say, degree $32$. We can give a heuristic argument as to why this is to be expected as follows. For this range of degrees, the majority of the vertex-transitive graphs are Cayley graphs, and so groups that contribute large numbers of Cayley graphs dominate the enumeration. If a group $G$ has $a$ involutions and $b$ non-involutions, then a first approximation to the number of distinct Cayley sets (i.e. inverse-closed subsets of $G$ that are pairwise inequivalent under $\rm{Aut}(G)$) is given by
\[
\frac{2^{a+b/2}} {| \rm{Aut}(G) |}.
\]
For a fixed degree, this value will be greatest when $a$ is large and $|\rm{Aut}(G) |$ is small. For the range of degrees currently under consideration, the \emph{dihedral groups} with relatively large $a$ and relatively small $b$ give the greatest value. If the degree is a power of $2$ then the elementary abelian $2$-group has the largest possible value of $a$, because every non-identity element is an involution. However the automorphism group of $\mathbb{Z}_2^n$ is $GL(n,2)$, which is very large, and so elementary abelian groups contribute relatively few Cayley graphs to the total.


\begin{table}
\begin{tabular}{|lrr|lrr|lrr|}
\hline
$n$&$t(n)$&$c(n)$&$n$&$t(n)$&$c(n)$&$n$&$t(n)$&$c(n)$\\
\hline
&&&$16$&$286$&$278$&$32$&$677402$&$659232$\\
&&&$17$&$36$&$36$&$33$&$6768$&$6768$\\
$2$&$2$&$2$&$18$&$380$&$376$&$34$&$132580$&$131660$\\
$3$&$2$&$2$&$19$&$60$&$60$&$35$&$11150$&$11144$\\
\hline
$4$&$4$&$4$&$20$&$1214$&$1132$&$36$&$1963202$&$1959040$\\
$5$&$3$&$3$&$21$&$240$&$240$&$37$&$14602$&$14602$\\
$6$&$8$&$8$&$22$&$816$&$816$&$38$&$814216$&$814216$\\
$7$&$4$&$4$&$23$&$188$&$188$&$39$&$48462$&$48462$\\
\hline
$8$&$14$&$14$&$24$&$15506$&$15394$&$40$&$13104170$&$13055904$\\
$9$&$9$&$9$&$25$&$464$&$464$&$41$&$52488$&$52488$\\
$10$&$22$&$20$&$26$&$4236$&$4104$&$42$&$9462226$&$9461984$\\
$11$&$8$&$8$&$27$&$1434$&$1434$&$43$&$99880$&$99880$\\
\hline
$12$&$74$&$74$&$28$&$25850$&$25784$&$44$&$39134640$&$39134544$\\
$13$&$14$&$14$&$29$&$1182$&$1182$&$45$&$399420$&$399126$\\
$14$&$56$&$56$&$30$&$46308$&$45184$&$46$&$34333800$&$34333800$\\
$15$&$48$&$44$&$31$&$2192$&$2192$&$47$&$364724$&$364724$\\
\hline
\end{tabular}
\vspace{2mm}
\caption{Numbers $t(n)$, $c(n)$ of transitive and Cayley graphs of order $n$}
\label{tab:transgraphs}
\end{table}

\section{Two Applications}

One of the major reasons to construct catalogues of combinatorial objects is to gather evidence relating to conjectures or other open questions. Even if a newly-constructed catalogue does not directly contain a counterexample to a conjecture (thereby immediately resolving it), it can be useful in refining a researcher's intuition regarding both the typical and extremal behaviour of the objects in the catalogue. 

In the remainder of this section, we consider two areas in which computational evidence has played a role, and augment that evidence with information derived from the list of transitive groups described in this paper.

\subsection{Elusive Groups}

A permutation group $G$ is called {\em elusive} if it contains no \emph{fixed-point-free} elements (i.e., \emph{derangements}) of prime order.  Elusive groups are interesting because of their connection to Maru\v{s}i\v{c}'s \emph{Polycirculant Conjecture} which asserts that the automorphism group of a vertex-transitive digraph is \emph{never elusive}.  In principle, a positive resolution of the polycirculant conjecture may simplify the construction and analysis of vertex-transitive graphs and digraphs, as it would then be possible to assume the presence of an automorphism with $n/p$ cycles of length $p$ for some prime $p$. For example, early catalogues of vertex-transitive graphs relied on ad hoc arguments to show that all transitive groups of the specific degrees under consideration have a suitable derangement of prime order.

A permutation group $G$ is called \emph{$2$-closed} if there is no group properly containing $G$ with the same orbitals as $G$. The automorphism group of a vertex-transitive digraph is necessarily $2$-closed, because it is \emph{already} the maximal group (by inclusion) that fixes the set of arcs of the digraph, which is a union of some of the orbitals. The conjecture can thus be strengthened to the assertion that there are no elusive $2$-closed transitive groups, yielding a conjecture that was first proposed by Klin.

One might hope  that there are simply no elusive groups at all, in which case both conjectures would hold vacuously, but in fact there are a number of sporadic examples of elusive groups and a handful of infinite families. All the known elusive groups are not $2$-closed, so do not provide counterexamples for either conjecture.  

It is relatively easy to test the groups for the property of being elusive by checking to see if any of the conjugacy class representatives are derangements of prime order. For the larger groups, it is often faster to first generate some number of randomly selected elements inside each of the Sylow subgroups in the hope of stumbling on a suitable derangement without the cost of computing all the conjugacy classes.

\begin{figure}
\begin{center}
\begin{tikzpicture}
\node (s1) at (0,4) {\small $((C_9^2 : Q_8) : C_3) : C_2$};
\node (s2) at (-3,2) {\small $(C_9^2 : Q_8) : C_3$};
\node (s3) at (3,2) {\small $(C_9^2 : C_8) : C_3$};
\node (s4) at (0,2) {\small $(C_9^2 : C_8) : C_2$};
\node (s2s) at (-2,0) {\small $C_9^2 : Q_8$};
\node (s3s) at (2,0) {\small $C_9^2 : C_8$};
\draw (s1)--(s2);
\draw (s1)--(s3);
\draw (s1)--(s4);
\draw (s2)--(s2s);
\draw (s3)--(s3s);
\draw (s4)--(s2s);
\draw (s4)--(s3s);
\end{tikzpicture}
\end{center}
\caption{The elusive groups of degree $36$}
\label{fig:elusive36}
\end{figure}
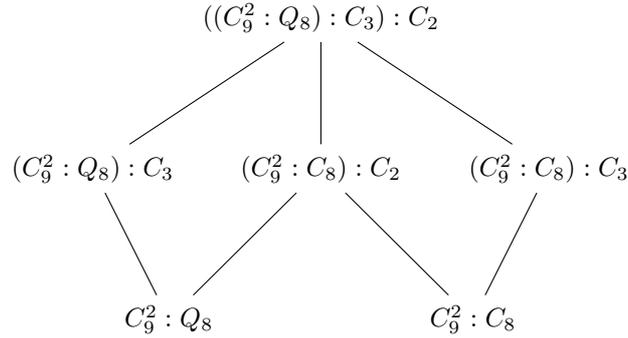

\begin{prop}
A transitive permutation group $G$ of degree $n < 48$ is elusive if and only if one of the following holds:
\begin{enumerate}
\item $G$ has degree $12$ and is either the group $M_{11}$ (acting on $12$ points), or one of its four proper transitive subgroups,
\item $G$ has degree $24$ and is one of the $19$ groups described by Giudici \cite{MR2356445},
\item $G$ has degree $36$ and is either a particular group of shape $((C_9^2 : Q_8) : C_3) : C_2$, or one of its five proper transitive subgroups, as shown in Figure~\ref{fig:elusive36}.
\end{enumerate}
\end{prop}

The complete lists of elusive groups for $n \leqslant 30$ were previously known (Giudici \cite{MR2356445}), and while some of the examples
for degree $36$ were known, the list was not complete. Prior to this work, the smallest degree for which 
the existence of an elusive group was undecided was $n=40$, which has now been ruled out.
The Appendix contains the indices of the elusive groups in the lists of transitive groups available in {\sc Magma}. 

\subsection{CI groups}

Recall that a \emph{Cayley digraph} for the group $G$ with connection set $C \subseteq G$ is the graph $\text{Cay}(G,C)$
with vertex set $G$, and where $(x, xc)$ is an arc of $G$ for all $c \in C$. Clearly $G$ acts regularly on $\text{Cay}(G,C)$ 
by left-multiplication, and any digraph admitting a regular automorphism group $G$ is a Cayley digraph for $G$. A Cayley digraph is loopless if and only if $\text{id}_G \not\in C$ and undirected if and only if $C = C^{-1}$. If $\sigma$ is an automorphism of the group $G$, then it is immediate that $\text{Cay}(G, C)$ is isomorphic to $\text{Cay}(G, C^\sigma)$. For some groups, these isomorphisms are the \emph{only} isomorphisms between the Cayley graphs, or digraphs,  on that group, in which case the group is called a CI-group or DCI-group respectively. (The acronym CI stands for ``Cayley Isomorphic'', and D for ``directed'' --- previously the terminology $\mathcal{G}$-CI group was used for the undirected case, but this seems to have fallen out of style.) We shall use \cite{MR1829620} for background.

One method of identifying the transitive graphs that are \emph{not} Cayley graphs would be to separately process each graph by computing its automorphism group and searching for regular subgroups. However it is considerably easier to precompute the entire collection of Cayley graphs for each of the groups of order less than $48$. Then the non-Cayley graphs are simply those that appear in the list of transitive graphs but not in the list of Cayley graphs.

As a side-effect, this process also allowed us to determine the entire collection of CI-groups of order up to $47$ in the following fashion. For each group $G$, we computed its automorphism group $A = \text{Aut}(G)$ and then, because we were primarily interested in the undirected case, determined the action of $A$ on the set
\begin{equation}\label{sdef}
S = \{ \{g, g^{-1} \} : g \in G, \ g \not= \text{id}_G\}.
\end{equation}
We used a straightforward orderly algorithm \cite{MR1614301} to calculate one representative from each orbit of $A$ on the power set $\mathcal{P}(S)$, thereby finding all the connection sets that are both inverse-closed and pairwise inequivalent under $A$. Finally, we constructed all the Cayley graphs with these pairwise-inequivalent connection sets and checked the lists for isomorphic pairs of graphs. The group $G$ is a CI-group if and only if this final step finds no isomorphic pairs of graphs. (If we were to repeat this for DCI-graphs, then $S$ would need to be the set containing all the non-identity elements of $G$, and the final step would construct all the Cayley digraphs.)

The result of this is a long list of CI and non-CI groups. However this can be presented in a more compact format by noting that a \emph{subgroup} of a CI-group is a CI-group and a \emph{quotient} of a CI-group is a CI-group (Dobson and Morris \cite{MR3338017}), thus permitting us to identify the \emph{minimal} non-CI-groups of order up to $47$.  In particular the following list of groups is the complete list of minimal non-CI groups of order at most $47$. 

\begin{enumerate}
\item The cyclic groups $C_{16}$, $C_{24}$, $C_{25}$, $C_{27}$, $C_{36}$, $C_{40}$, $C_{45}$.
\item The dihedral groups $D_{12}$, $D_{16}$, $D_{20}$, $D_{28}$, $D_{44}$.
\item One of $16$ groups listed individually in Table~\ref{noncigroups}.
\end{enumerate}

\begin{table}
\begin{tabular}{|ccc|ccc|}
\hline
Order & No. & Structure & Order & No. & Structure\\
\hline
$8$&$2$&$C_4\times C_2$&
$16$&$8$&$\mathit{QD}_{16}$\\
$16$&$9$&$Q_{16}$&
$18$&$3$&$C_3\times S_3$\\
$20$&$3$&$C_5:C_4$&
$24$&$3$&$\mathsf{SL}(2,3)$\\
$24$&$10$&$C_3\times D_8$&
$24$&$12$&$S_4$\\
$24$&$13$&$C_2\times A_4$&
$27$&$2$&$C_9\times C_3$\\
$27$&$3$&$(C_3\times C_3):C_3$&
$27$&$4$&$C_9:C_3$\\
$36$&$9$&$(C_3\times C_3):C_4$&
$36$&$11$&$C_3\times A_4$\\
$40$&$10$&$C_5\times D_8$&
$42$&$1$&$(C_7:C_3):C_2$\\
\hline
\end{tabular}
\vspace{2mm}
\caption{Minimal non-CI groups with $|G| < 48$ (not cyclic or dihedral).}\label{noncigroups}
\end{table}



The column labelled ``Structure'' is simply the output of GAP's {\tt StructureDescription} command. As this is not unique, we also include (in the column labelled ``No.'')  the number of the group in the list of \emph{small groups} of that particular order. These lists of small groups are found in both {\sc Magma} and GAP and, at least for the orders we are considering here, the numbering is consistent between the two programs. The \emph{regular representations} of these small groups are transitive groups of course, and so they appear in the lists of transitive groups of each degree. However the ordering of the regular transitive groups is \textit{\textbf{not the same}} as the ordering of the small groups. Therefore {\tt SmallGroup(deg,k)} and {\tt TransitiveGroup(deg,k)} are usually not isomorphic.

%

\subsection*{Acknowledgements}

We would like to thank Michael Giudici for helpful discussions regarding elusive groups and the polycirculant conjecture, and Joy Morris for helpful discussions regarding CI-groups.

\bibliographystyle{plain}
\bibliography{../../gordonmaster}

\appendix

\section{CI and non-CI groups}

In this section, we record in more detail some of the computations underlying the list of CI and non-CI groups.

Recall that a group $G$ is a CI-group if the only isomorphisms between Cayley graphs for $G$ are those induced by automorphisms of $G$ acting on possible connection sets (which in our case are inverse-closed sets without the identity).

Letting the term \emph{Cayley set} refer interchangeably to an ${\rm Aut}(G)$-orbit of subsets of $V(G)$ or an arbitrary member of such an orbit, we see that a group is CI if and if it has an equal number of Cayley sets and Cayley graphs. 

The number of Cayley sets of a group $G$ can be calculated using the \emph{cycle index polynomial} of $A = \rm{Aut}(G)$ acting on the set $S$ (defined above in \eqref{sdef}).  
This polynomial is the multivariate polynomial in variables $\{X_1, X_2, \ldots \}$
given by
\[
Z(A;X_1, X_2, \ldots) := \frac{1}{|A|} \sum_{a \in A} \prod X_k^{n_k(a)}.
\]
where $n_k(a)$ is the number of cycles of length $k$ in the permutation $a$.   The number of Cayley sets is then just an evaluation of this polynomial, in particular the evaluation $Z(A;2, 2, \ldots)$. The {\sc Magma} function below performs these computations using the (undocumented) function {\tt CycleIndexPolynomial}.

\bigskip

\hrule
\bigskip

{\small 
\begin{minted}{matlab}
numCayleySets := function(g)

  ptsinv := Setseq({ {x,x^-1} : x in g | Order(x) eq 2 });
  ptsnot := Setseq({ {x,x^-1} : x in g | Order(x) gt 2 });

  pts := ptsinv cat ptsnot;  generators := [];

  sym := Sym(#pts);

  autg := AutomorphismGroup(g);

  generators := [];
  for gen in Generators(autg) do
    Append(~generators,[Position(pts,gen(pts[i])) : i in [1..#pts]]);
  end for;

  autonpts := sub<sym|generators>;
  ci := CycleIndexPolynomial(autonpts);

  return Evaluate(ci,[2 : i in [1..Rank(Parent(ci))]]);

end function;
\end{minted}
}

\bigskip
\hrule

By further refining the cycle index polynomial to count separately the numbers of cycles of each length that consist of elements of $S$ corresponding to involutions, and elements of $S$ corresponding to a non-involution along with its inverse, it would be possible to count the number of Cayley sets of each valency. This would allow some greater exploration of the more refined property of being $m$-CI, which means that the number of Cayley sets and Cayley graphs of degree $m$ coincide.

\end{document}